\documentclass[12pt]{article}
\usepackage{setspace}
\usepackage{textcomp}
\usepackage{amsmath,amsthm,amssymb}
\usepackage[left=1.5in,top=1in,right=1in,nohead]{geometry}
\usepackage{latexsym}
\usepackage{mathrsfs}
\usepackage{rawfonts}
\usepackage[dvips]{graphicx}
\usepackage{stmaryrd}

\usepackage[sc,osf]{mathpazo}

\def\bct{\begin{center}}
\def\ect{\end{center}}
\def\beg{\begin}

\def\<{\langle}
\def\>{\rangle}

\def\mbb{\mathbb}
\def\mbbz{\mathbb Z}

\def\ni{\noindent}

\def\tn{\textnormal}

\title{Locally conformally flat and self-dual structures on simple 4-manifolds} \author{Mustafa Kalafat}
\begin{document}
\maketitle
\begin{abstract} This is a survey article on the existence of locally conformally flat (LCF) and self-dual (SD) metrics on various basic 4-manifolds like simply-connected ones or product types.
\end{abstract}
%\keywords{Locally conformally flat manifolds, self-dual manifolds, K
%\" ahler manifolds}

%realizing rich topological types, which are obtained
%from conformal compactification of the 3-manifolds, that are built
%from the Panelled Web Groups. These manifolds have strictly negative
%scalar curvature

%\tableofcontents

\section{Introduction}

\vspace{.05in}
A Riemannian manifold $(M,g)$ is called {\em locally conformally flat (LCF)} if there is a neighborhood $U$ of any point, and a strictly positive smooth function $f$ such that $\tilde{g}=fg$ is a metric of zero sectional curvature everywhere on $U\subset M$. %\vspace{.05in}
Sometimes the more concise terminology {\em conformally flat} is also used since globally conformally flat manifolds admit flat metrics. They are quotients of $\mathbb R^n$ so that there is no need to give them a new name. 
Alternating and symmetry properties imply that we can think of the $(0,4)$ 
Riemann curvature tensor $R$ as an element of the space $S^2\Lambda^2M \subset \otimes^4T^*M$. It also satisfies the algebraic Bianchi identity, hence it lies in the kernel of the  {\em Bianchi symmetrization map} 
$$b: S^2\Lambda^2M \to S^2\Lambda^2M,  \hspace{6mm} b(T)(x,y,z,t):={1 \over 3}\,T((x,y,z),t).$$
Since $b^2=b$ and $b$ is ${GL}(T^*M)$-equivariant, we have the equivariant decomposition $S^2\Lambda^2M=\tn{Ker}\,b\oplus\tn{Im}\,b$, where we call this kernel  as the space of {\em curvature-like tensors}. Thinking $\tn{Ker}\,b$ as an invariant $O(g)$-module, we have a unique irreducible decomposition, according to this the curvature tensor decomposes as $R=U\oplus Z \oplus W$. The components can be computed as 
$$U={s\over 2n(n-1)}g\varowedge g \hspace{.5cm}\textnormal{and} \hspace{.5cm} Z={1\over n-2}\stackrel{\circ}{\tn{Ric}}\varowedge \, g   $$ 
where $s$ is
the scalar curvature, $\stackrel{\circ}{\tn{Ric}}\,=\tn{Ric}-{s\over n}g$ is the
trace-free Ricci tensor, ``$\varowedge$" is the {\em Kulkarni-Nomizu product} defined by,
$$A\varowedge B\,(X,Y,Z,T):=
\beg{array}{|cc|} A(X,Z) & B(X,T)\\ A(Y,Z)& B(Y,T) \end{array} \,+\, 
\beg{array}{|cc|} B(X,Z) & A(X,T)\\ B(Y,Z)& A(Y,T) \end{array} \, , $$
which is commutative and multiplies two symmetric 2-tensors to produce a curvature-like 4-tensor. Finally $W$ is the {\em Weyl tensor} which is defined to be what is left over from the first two pieces.

Next, assume that we are in dimension $n=4$ and our manifold is oriented. Then the metric together with the orientation determines a unique volume form $\omega_g$. In this case we define the {\em Hodge star involution}  $*_g: \Lambda^2\to\Lambda^2$ pointwise by imposing the equality $\langle \alpha,\beta \rangle \omega_g=\alpha\wedge*\beta$. 
This yields the $\pm 1$ eigenspace decomposition $\Lambda^2=\Lambda^2_+ \oplus \Lambda^2_-$ of the $2$-forms. These $3$-dimensional eigenspaces are interchanged if one works with the reversed orientation. 
By an appropriate change of indices, we consider $W:\Lambda^2\to\Lambda^2$ as an operator. One can show that the mixed parts $W_+^- : \Lambda^2_+ \to \Lambda^2_-$ and $W_-^+ : \Lambda^2_- \to \Lambda^2_+$ vanish so that the Weyl tensor decomposes as $W=W_+^+\oplus W_-^-$. Abbreviating $W_\pm=W_\pm^\pm$, we say that the Riemannian manifold is {\em self-dual} if $W_-\equiv 0$, {\em anti-self-dual} if $W_+\equiv 0$ respectively. We also call each of these two cases as {\em half-conformally flat} if we do not want to specify any orientation. This terminology actually comes from the interpretation that $W$ is a conformally invariant tensor if one considers it as a $(1,3)$ tensor. One can show that the manifold is conformally flat if and only if the tensor $W\equiv 0$ for dimensions $n\geq 4$.  See \cite{kuhnel} for a proof. In dimension $n=3$ local conformal flatness is determined by the Schouten tensor, and in dimension $n=2$ all manifolds are locally conformally flat. Basic examples of LCF manifolds are constant sectional curvature spaces, e.g., $S^n$ and $T^n$, their products with $S^1$ or 
$\mathbb R$, and products of two Riemannian manifolds with constant sectional curvature $1$ and $-1$ respectively. See \cite{besse} for further details. 
Note that SD metrics are only defined on orientable $4$-manifolds. 
Basic SD $4$-manifolds are first of all LCF spaces, e.g. $S^4$, $T^4$, $S^1\times S^3$ (see Theorem \ref{S1xS3}), among non-LCF ones we have $\mathbb{CP}_2$ with its Fubini-Study metric (see Corollary \ref{cp2fubinistudy}) and scalar-flat-K\"ahler (SFK) surfaces (see Theorem \ref{sfk=>asd}) in particular the K3 surface. SD $4$-manifolds are first considered by \cite{penrose} and later put on a firm mathematical foundation in \cite{ahs}.

\vspace{.05in} %Simdi  SD tanimi. $S^2\times S^2,S^1\times S^3$ Besse - AHSden tiyolar.

In this survey, we analyze LCF and SD structures on various simple $4$-manifolds like product type or simply-connected. These results are spread out various places. Some of them are not written, if so not in detail. We hope that it is a good public service to accumulate these results in an article. 
Interested reader may consult to the resources \cite{handle} and \cite{sd} for some 
recent progress on this type of geometry. 
In section \ref{sectools} we introduce the basic tools, 
in section \ref{secsimpleieproductanssimplyconnected} we apply these 
tools on the manifolds and finally in the appendix we present a partial result along with an open problem. 

\vspace{.05in}

{\bf Acknowledgements.} Thanks to S. Finashin, C. Koca and M. Korkmaz for useful discussions.

\section{Tools}\label{sectools}

In this section we develop the main tools to analyze our spaces. 
Let $M$ be a closed, oriented $4$-manifold. We have the following two basic identities which connect quantities related to curvature with topological numbers. 

\beg{equation}\label{euler} \chi(M)={1\over 8\pi^2}\int_M {s^2\over 24}- %{ |\circ{\tn{Ric}}|\over 2} 
{|\stackrel{\circ}{\tn{Ric}}|^2 \over 2}+ |W_+|^2+|W_-|^2 \omega_g \end{equation}

\beg{equation}\label{signature} \tau(M)={1\over 12\pi^2}\int_M |W_+|^2-|W_-|^2 \omega_g \end{equation}

\ni The first one is called the 
{\em generalized Gauss-Bonnet theorem} \cite{allendoerferweil,singerandthorpe} which can also be generalized to all even dimensions.  % Hirzebruch top meth Referanslari dogru veren:  seminare bourb claude
The second one is called the {\em signature formula} which is specific to dimension $4$. It is obtained through the {Hirzebruch signature theorem} \cite{hirzebruchsignaturetheorem} 
$\tau(M)={1 \over 3}p_1[M]$ by expanding the Pontrjagin class with Chern-Weil theory \cite{singerandthorpe}. See also second volume of \cite{kn}. 
These two equations relate geometric information with the topological ones. As an immediate application for example if a $4$-manifold admits a locally conformally flat metric then both of the self-dual and anti-self-dual Weyl curvatures vanish since $W=W_+\oplus W_-=0$. Hence the signature formula (\ref{signature}) implies that the signature $\tau(M)=0$. Perhaps this signature condition is the most important topological obstruction for admitting LCF metrics. 
We start with a very common way of producing ASD metrics on a K\"ahler manifold. Recall that a Riemannian manifold is called {\em scalar-flat} 
if its scalar curvature is zero everywhere.

%\beg{thm}\label{sfk=>asd} Scalar-flat K\"ahler(SFK) surfaces are anti-self-dual(ASD).
%\end{thm}

\beg{thm}\label{sfk=>asd} A K\"ahler surface is scalar-flat(SF) iff anti-self-dual(ASD).
\end{thm}
\beg{proof} We follow \cite{besse} and \cite{clexdef}. The alternating and symmetry properties of the curvature tensor implies that it is a  symmetric section of the bundle $\wedge^2 \otimes \wedge^2$. Furthermore, since we are on a K\"ahler manifold we have the identities \cite{kn}, 
$$R(JX,JY)Z=R(X,Y)Z \hspace{3mm} \tn{and} \hspace{3mm} R(X,Y)JZ=JR(X,Y)Z.$$
These imply respectively the $J$-invariance of the first and the second pair of components of the curvature tensor hence it is a type $(1,1)$ real $2$-form in both of these components. So we can think of the K\"ahler curvature tensor as a symmetric section of $\wedge^{1,1} \otimes \wedge^{1,1}$ or after dualizing the second component, as a symmetric element of $\tn{End}(\wedge^{1,1})$.
%We follow \cite{clexdef}. 

Since we are on a complex manifold, we have the Dolbeault decomposition of $2$-forms $\wedge^2_{\mbb C}=\wedge^{2,0}\oplus \wedge^{1,1} \oplus \wedge^{0,2}$ according to their type. This decomposition is orthogonal with respect to the action of the Hodge star operator $*_g$. We also have the orthogonal eigenspace decomposition in dimension $4$ as explained in the introduction. In the K\"ahler case, if one complexifies these eigenspaces, we claim to have the following:
%Here, the subspaces $\wedge^{1,1}$ and 
%$\wedge^{2,0} \oplus \wedge^{0,2}$ are invariant under the Hodge star operator $*_g$.
$$\wedge^2_{+\mbb C} = \mbb C \omega \oplus \wedge^{2,0}\oplus \wedge^{0,2}$$
$$\wedge^2_{-\mbb C} = \wedge^{1,1}_0$$
Here, the set of {\em primitive $(1,1)$ forms} $\wedge^{1,1}_0$ can be defined to be the orhogonal complement of $\omega$ in $\wedge^{1,1}$. 
Take a unitary coframe $\{dz^1,dz^2\}$ at a point. Then the K\"ahler form is
$$\omega={i \over 2}\sum_{i=1}^2 dz^i\wedge d\bar{z}^i = dx^1\wedge dy^1+dx^2\wedge dy^2$$
and the volume form is computed as
$$\omega_g={\omega^2 \over 2!}
={ -1 \over 4}\, dz^{1\bar 12 \bar 2} 
=dx^1\wedge dy^1\wedge dx^2\wedge dy^2.$$ 
Since $*_g(dx^i\wedge dy^i)=dx^{3-i}\wedge dy^{3-i}$ we have $*_g \omega=\omega$, i.e., the K\"ahler form is self-dual. 
This observation links to the following interpretation of the primitive $(1,1)$ forms. 
We can write $\wedge^{1,1}_0=\tn{Ker}\,L$,
for the {\em Lefschetz operator} $L:\wedge^{1,1}\to\wedge^{2,2}$ defined by $L(\alpha)=\alpha\wedge\omega$ since
$$\langle\alpha,\omega\rangle=0 \Leftrightarrow 0=\langle\alpha,\omega\rangle \omega_g=\alpha\wedge *_g \omega=\alpha\wedge\omega=L(\alpha).$$
The forms $\omega_5=dz^1\wedge dz^2$ and $\omega_6=d\bar{z}^1\wedge d\bar{z}^2$ pointwise generate the complex bundles $\wedge^{2,0}$ and $\wedge^{0,2}$ respectively. One easily checks that their real parts are self-dual.
Hence  $\omega$, $\omega_5$, $\omega_6\in \wedge^2_{+\mbb C}$ are orthogonal and their real parts are elements of the real rank $3$ vector bundle of self-dual $2$-forms $\wedge^2_+$. On the other hand the $(1,1)$ forms $\omega_2=dz^1\wedge d\bar z^2$,  $\omega_3=dz^2\wedge d\bar{z}^1$ and 
 $\omega_4=dz^1\wedge d\bar z^1-dz^2\wedge d\bar z^2$ are orthogonal and their real parts are anti-self-dual. For example 
 $*_g\Re\omega_2=*_g(dx^{12}+dy^{12})=-(dx^{12}+dy^{12})=-\Re\omega_2$. Moreover they are all orthogonal to 
$\omega$, $\omega_5$, $\omega_6\in \wedge^2_{+\mbb C}$. %$\wedge^2_{+\mbb C}$.
 
 Now, since the curvature operator $\mathcal R$ is in $\tn{End}(\wedge^{1,1})$, its upper left piece $\mathcal R_+^+=\mathcal W_+ + {s \over 12} I$ is an element of $\tn{End}(\mathbb C \omega)$. So suppose $\mathcal R_+^+=f \, \omega\otimes\omega^\sharp$ for some function $f: M\to \mbb R$.

{ \renewcommand*{\arraystretch}{1.4}
 \begin{table}[!h]
\centering
$\beg{array}{|c|c|c|} \hline 
\mathcal R& \mbb C\omega \oplus \wedge^{2,0} \oplus \wedge^{0,2} & \wedge^{1,1}_0 \\ %[1ex]
\hline
\mbb C\omega \oplus \wedge^{2,0} \oplus \wedge^{0,2} & \mathcal W_+ + {s \over 12} I & \stackrel{\circ}{\tn{Ric}}\\
\hline
\wedge^{1,1}_0 & \stackrel{\circ}{\tn{Ric}} &  \hspace{.5cm} \mathcal W_- + {s \over 12} I \hspace{.5cm} \\
\hline
\end{array}$ 
\vspace{2mm} 
\caption{Curvature operator for K\"ahler surfaces}
\label{tab:curvkahler}
\end{table}
}

{ \renewcommand*{\arraystretch}{1.6}  
\begin{table}[h]
\centering
$\beg{array}{|c|c|@{}c@{}|} \hline
\mathcal R & \mbb C\omega & \wedge^{1,1}_0 \\
 \hline
\mbb C\omega  &  s/4\cdot & \stackrel{\circ}{\tn{Ric}} \\
 \hline
\wedge^{1,1}_0 & \stackrel{\circ}{\tn{Ric}} & \;\; \mathcal W_- + {s \over 12} I \;\; \\
 \hline
\end{array}$ \hspace{1.5cm} 
$\beg{array}{|c|c|c|c|} \hline
\mathcal R_+^+& \mbb C\omega & \wedge^{2,0} & \wedge^{0,2}\\
 \hline
\mbb C\omega  &  s/4\cdot & 0 & 0 \\
 \hline
\wedge^{2,0} & 0 & 0 & 0 \\
 \hline
\wedge^{0,2} & 0 & 0 & 0 \\
 \hline
\end{array}$
\vspace{2mm}
\caption{Curvature operator and its self-dual part for K\"ahler surfaces}
\label{tab:sdcurvkahler}
\end{table}
 }

\noindent To figure out the function $f$ we need to compute some inner products. First compute the norm of the K\"ahler form $\omega$. Since the volume form $\omega_g=\omega^2/2!$\, we have $$\langle\omega,\omega\rangle\,\omega_g = \omega\wedge *\omega=\omega\wedge\omega=2\,\omega_g.$$
Hence we get $|\omega|=\sqrt{2}$.
\footnote{In general dimension $n$, writing in local orthonormal frame $\omega=e_1\wedge e_2+...+e_{2n-1}\wedge e_{2n}$ we compute $|\omega|=\sqrt{n}$. Then from 
$|\omega|^2\omega_g=\omega\wedge*\omega$ we get $*\omega=|\omega|^2\omega^{n-1}/n!$  
so $*\omega=\omega^{n-1}/(n-1)!$} 
Secondly we want to compute the inner product 
$\langle\rho,\omega\rangle$\, where $\rho$ is the {\em Ricci form} defined by 
$\rho(\cdot,\cdot)=\tn{Ric}(J\cdot,\cdot)$. The above trick does not work in this case. We need to use some tensor analysis. Start with fixing a convention for the complex structure. 
Suppose $J=J_i^{\phantom{k}k}\,  dx^i \otimes \partial_k$. Then the basic property $J^2=-Id$ reads $J_i^{\phantom{i}k}J_k^{\phantom{j}j}=-\delta_i^j$ 
\footnote{If one fixes the alternative convention $J=J^k_{\phantom{i}i}\, \partial_k\otimes dx^i$ one gets $J^j_{\phantom{i}k}J^k_{\phantom{i}i}=-\delta_i^j$ and $\omega_{ij}=-J_{ij}$, a negative sign. } 
in terms of tensors. Keep %ing 
in mind that the $\omega$ and $J$ %have the same indices since they correspond to one another under the metric and they 
are skew symmetric tensors. We compute the following:
$$\omega^{ij}=\omega(dx^i,dx^j)=g(Jdx^i,dx^j)=g(J_k^{\phantom{i}i} dx^k,dx^j)
=J_k^{\phantom{i}i} g^{kj}=J^{ji}=-J^{ij}.$$
$$\rho_{ij}=\rho(\partial_i,\partial_j)=\tn{Ric}(J\partial_i,\partial_j)
=J_i^{\phantom{k}k}\tn{Ric}(\partial_k,\partial_j)=J_i^{\phantom{k}k} R_{kj}.$$
And then,
\begin{align*} 
\<\rho,\omega \>&={1 \over 2!}\, \rho_{ij}\, \omega^{ij}
={ 1 \over 2}\,J_i^{\phantom{k}k} R_{kj} (-J^{ij})
={-1 \over 2}\,J_i^{\phantom{k}k} R_{kj} J_k^{\phantom{j}j} g^{ki} \\
&={-1 \over 2}\, (-\delta_i^j) R_{kj} g^{ki}
={ 1 \over 2}\, R_{ki} g^{ki}
={ s \over 2}.
\end{align*}
 
Now, writing $\mathcal R \omega = f \omega + g \omega^\perp$, multiplying both sides with the K\"ahler form and using $\mathcal R \omega = \rho$, we get the following:
$$\beg{array}{rcl}
\langle \mathcal R \omega , \omega \rangle &=& f\, \langle \omega , \omega \rangle \\
\langle \rho , \omega \rangle &=& f\, \langle \omega , \omega \rangle\\ 
{s/ 2} &=& f\cdot 2\\
{s/ 4} &=& f.     \end{array}$$
Hence $\mathcal W_+$ is a multiple of the scalar curvature $s$. Therefore we have $s=0$ if and only if $\mathcal W_+=0$. \end{proof}

This result has an immediate corollary.

%\cite{} citationu nedir bunun, Berger veya Bergery felan mi? LSde bulamadim, Besse veya baska Claude lazim.
\beg{cor}\label{kahlerWplusandscalarnorms} For K\"ahler metrics on a complex surface we have the pointwise identity   
\beg{equation} %%%\label{kahler}
 |W_+|^2={s^2 \over 24}.  \label{eq:kahlerWplus} \end{equation}
\end{cor}

\beg{proof} Since in the K\"ahler case the term $\mathcal R_+^+=W_+ + {s \over 12} I$ acts by the following
    $${s \over 4}\left[\beg{array}{ccc}1&0&0\\ 0&0&0\\ 0&0&0 \end{array}\right],$$ 
as in the Table \ref{tab:sdcurvkahler}, we have $$W_+={s \over 12}\left[\beg{array}{ccc}2&0&0\\ 0&-1&0\\ 0&0&-1 \end{array}\right].$$ Taking the norms of both sides yields the result.
\end{proof}

\beg{thm}\label{KE} If a $4$-manifold $M$ admits a K\"ahler-Einstein (KE) metric $g$ then \linebreak
$\chi(M)=3\tau(M)$ if and only if $g$ is self-dual (SD).  %--ASD-- olmuyor duz orientasyonda self-dual iste 
\end{thm}

\beg{proof} For K\"ahler metrics we have the pointwise identity   \eqref{eq:kahlerWplus}
\beg{equation}\label{kahler} |W_+|^2={s^2 \over 24}. \nonumber \end{equation}
For Einstein metrics, by definition the trace-free Ricci tensor $\stackrel{\circ}{\tn{Ric}}\, =0$.
Plugging these two identities into the Gauss-Bonnet formula (\ref{euler}) we get
\beg{equation}\label{eulerKE} 8\pi^2\chi=2\,\|W_+\|^2+\|W_-\|^2.\end{equation}
Eliminating $\|W_+\|$ from (\ref{signature}) and (\ref{eulerKE}) we get the following equality  $$8\pi^2(\chi-3\tau)=3\,\|W_-\|,$$ which yields the result. \end{proof}

Finally we state the celebrated theorems of Kuiper. See also \cite{howard} for a recent exposition and improvement.

\beg{thm}[\cite{kuiper}] \label{kuiperthm} Let $(M^n,g)$ be a simply connected, LCF $n$-manifold of class $C^1$. Then there is a conformal immersion $f:M\to S^n$. If in addition $M$ is compact, then this map is a conformal diffeomorphism.
\end{thm}

Here comes another very useful theorem of Kuiper. See also \cite{osamukobayashi}.
 
\beg{thm}[\cite{kuiper2}] \label{kuiperthm2} Universal cover of a compact, LCF space with an infinite Abelian fundamental group must be $\mbb R^n$ or $\mbb R \times S^{n-1}$.
\end{thm}

\section{Simple 4-manifolds}\label{secsimpleieproductanssimplyconnected}
%\section{Product 4-manifolds}\label{secsimpleieproductanssimplyconnected}

We start with a basic space, the complex projective space with its standart Fubini-Study metric. This can be though as the metric quotient of $S^{2n+1}\subset \mbb C^{n+1}$ by unit complex scalar multiplication. Alternatively on $\mbb C^n$ take the following complex metric coefficients:
$$g_{FS_{i\bar{j}}}=g_{FS}(\partial_i,\partial_{\bar{j}})
:={(1+|z|^2)\delta_{i\bar{j}}-z_i z_{\bar{j}} \over(1+|z|^2)^2}.$$
Taking the completion of this space gives the complex projective space 
$\mbb{CP}_n$. As an application of the theorems in the previous section we obtain the following.  
\beg{cor}\label{cp2fubinistudy} 
$(\mbb{CP}_2,g_{FS})$, the complex projective space with its Fubini-Study metric is self-dual.
\end{cor} 
\beg{proof}One can easily compute $\chi(\mbb{CP}_2)=3$. Since the intersection form 
$Q_{\mbb{CP}_2}=[1]$ we have $\tau=1$. These satisfy the equality in Theorem \ref{KE}.\end{proof}

\beg{thm} The underlying smooth manifold of $\mbb{CP}_2$, the complex projective space does not admit any LCF metrics.
\end{thm}
\beg{proof}The basic obstruction signature $\tau(\mbb{CP}_2)=1$ is nontrivial.\end{proof}

Next we work on the 4-manifold $S^2\times S^2$, the product of two spheres. 
The  spheres $S^2\times q$ and $p\times S^2$ generating the homology have self-intersection zero, and $+1$ with each other. So 
$$Q_{S^2\times S^2}=\left[ \begin{array}{cc} 0 & 1\\ 1 & 0\end{array} \right],$$
called the hyperbolic matrix and denoted by $H$. It has eigenvalues $\pm 1$ and hence
the signature $\tau=0$. 
\footnote{Alternatively, the map $(I_3,-I_3)$ is an orientation reversing diffeomorphism of $S^2\times S^2$. So 
Hirzebruch signatures mapped onto each other and $\tau=-\tau$.} 
Hence an obstruction vanishes for locally conformally flatness. However, this turns out to be not sufficient as follows.

\beg{thm}\label{spheretimessphere} The 4-manifold $S^2\times S^2$ does not admit any LCF nor even SD metrics. %(tauyla felan mi state etsek?)
\end{thm}
\beg{proof} Since by Kuiper's theorem \cite{kuiper}, any compact, simply-connected, LCF Riemannian 4-manifold is conformally equivalent to the round 4-sphere, 
$S^2\times S^2$ does not admit any LCF metric. 
Suppose it does have a self-dual metric. Then $W_-=0$ and by the signature formula (\ref{signature}) the integral $0=\int_M |W_+|^2\omega_g$. So, pointwise $W_-=0$ and $W=0$. 
This yields a LCF metric which is already a contradiction by the previous argument.\end{proof}

This example also illustrates the fact that the product of LCF manifolds may not be LCF.

%\ni {\bf Q:} Hangi blowupi ilk SD admit eder?

\beg{cor}\label{K3} The 4-manifold $K3$ is not LCF, but SD.\end{cor}
\beg{proof} Since the intersection form of a K3 surface is $Q=2E_8\oplus 3H$, the signature $\tau(K3)$ is $-16$ which is nonzero so that it can not be LCF. Since the first Chern class $2\pi c_1({K3})$ is zero, the zero form is a $(1,1)$-form representing this class. Since $K3$ is compact and K\"ahler 
by Yau's resolution to the Calabi's problem \cite{yau},
 there is a K\"ahler metric $g$ in the same class with Ricci form $\rho_g\equiv 0$. Hence $g$ is Ricci flat, so scalar flat. SFK surfaces are ASD, hence SD with the reversed orientation.  \end{proof}

\beg{thm} The smooth $4$-manifold $\mbb{CP}_2\# \,\overline{\mbb{CP}}_2$ does not admit any LCF nor SD metrics.
\end{thm}
\beg{proof} Since $\tau=0$, the argument in the proof of Theorem \ref{spheretimessphere} is again valid for this 4-manifold.
% o halde LCF mi? Burns?
\end{proof}
A similar argument also excludes the manifolds $k\mbb{CP}_2\# \,k\overline{\mbb{CP}}_2$ for any $k>0$. 
Since we used it multiple times, it is convenient to sum up the idea of the proof as follows.
\beg{thm}\label{signature0} Let $M$ be a compact, oriented 4-manifold  with signature $\tau=0$, then we have the following.

1. A metric $g$ on $M$ is SD iff LCF.

2. If $M$ is simply connected but not diffeomorphic to $S^4$ then it does not admit any LCF or SD metric.

\end{thm}

\beg{thm} The smooth $4$-manifold $S^1 \times S^1 \times S^1 \times S^1 = T^2\times T^2$ is LCF hence SD.
\end{thm}
\beg{proof} Since $T^4$ is a quotient of  $\mbb R^2$ by $\mbbz^2$ isometries, it is flat hence LCF and consequently SD. 

The second assertion is alternatively seen as follows: The product of flat metrics on the components is scalar-flat-K\"ahler (SFK). So by the Theorem \ref{sfk=>asd} it is ASD.\end{proof}

As a by-product of local conformal flatness we also compute $\tau(T^4)=0$. Bieberbach's theorem states that the only compact %Riemannian 
manifolds that admit flat metrics are $T^n$ and its finite quotients. Hence these are the only (globally) conformally flat manifolds. Since these are already classified, the phrase ``conformally flat" is usually used in place of LCF in the literature. 
%\ni {\bf Q:} {\em Globally conformally flat} compact $T^4$.
%Let us now finally check the last possible torus in dimension four.
Before checking another possible basic product in dimension four, let us look at its universal cover.

\beg{thm} The smooth $4$-manifold $S^2 \times \mbb R^2$ is LCF hence SD.
\end{thm} \beg{proof}%x Since both conditions are local, we can work on $S^2\times\mbb R^2$ instead. 
Since $S^2\times \mbb R$ is diffeomorphic to $\mbb R^3-\{0\}$, the inherited standart metric is of constant (positive) sectional curvature. A constant sectional curvature space times $\mbb R$ is LCF. \end{proof}
%Similarly this computes $\tau(S^2\times T^2)=0$ as well.

Next result indicates that the above manifold can not have isometries that give $S^2 \times T^2$ as the quotient.

\beg{thm} The smooth $4$-manifold $S^2 \times T^2$ is neither LCF nor SD.
\end{thm}
\beg{proof} Because of the antipodal orientation reversing isometry of the $2$-sphere, we have the signature $\tau(S^2 \times T^2)=0$. Hence this is LCF iff SD.
Secondly, since this manifold has Abelian infinite fundamental group $\mbbz^2$, if it admits a LCF metric then applying  Kuiper's second theorem (Theorem \ref{kuiperthm2}) it should have universal cover $\mbb R^4$ or $\mbb R\times S^3$. Since the universal cover is $S^2\times \mbb R^2$ this gives a contradiction.  In general this argument applies to $S^p\times T^q$ for $p,q \geq 2$. 
\end{proof}
This manifold is particularly interesting since even though it is not LCF, the following infimum of all the possible Weyl energies on the manifold,
$$W(M):=\tn{inf} \left\{ \int_M |W_g|_g^2 \,\omega_g : g\in \mathcal{M}_M \right\}$$
called the {\em Weyl invariant} is zero, where $W_g$ is the conformally invariant type $(1,3)$ Weyl tensor, $\mathcal{M}_M$ is the space of smooth metrics. See \cite{osamukobayashi} for a proof. Hence the infimum is not attained.

\vspace{.05in}

%\section{Other basic manifolds} %\vspace{.05in}
One can easily compute the Euler characteristic of $S^1\times S^3$ to be $\chi=0$, since the $S^3$ component is an oriented closed 3-manifold hence $\chi(S^3)=0$. 
Since there is no second homology, the signature $\tau=0$. 
In this case it is a candidate of a manifold which may admit LCF metrics. This turns out to be the case as follows. 

\beg{thm}\label{S1xS3} 
The 4-manifold $S^1\times S^3$ with its standard product metric is LCF and hence SD.\end{thm}

\beg{proof} The proof is adapted from \cite{jeffdg}. Locally we can think this as $\mbb R\times S^3$, where $\mbb R$ is the flat line with $g_{\mbb R}=dx\otimes dx$ and $S^3$ has the round metric of curvature $+1$. Rm denoting the Riemann curvature tensor, $\varowedge$ is the Kulkarni-Nomizu product which is commutative,
$$\renewcommand{\arraystretch}{2}
\beg{array}{rcl}
\tn{Rm}_{\mbb R \times S^3} & = & \tn{Rm}_{\mbb R} + \tn{Rm}_{S^3} \\ 
  & = & 0+{1\over 2}\, g_{S^3}\varowedge    % \wedge   \rlap{\hspace{.05cm}$\wedge$}D\rlap{$\wedge$}            {\tn{O}}  
  g_{S^3}\\
  &=&{1\over 2}\, (g_{S^3}+dx^2)\varowedge(g_{S^3}-dx^2)\\
  &=&{1\over 2}\,  g_{\mbb R \times S^3}\varowedge(g_{S^3}-dx^2)\\  
  &=& \Psi({1\over 2}\, (g_{S^3}-dx^2))  
\end{array}$$
where by the definition of the product $dx^2\varowedge dx^2=0$, and 
$\Psi : S^2(T^*M) \longrightarrow \tn{Ker}\,b$ is defined by $\Psi(h)=h \varowedge g_M$. Here, $\tn{Ker}\,b$ is the space of curvature-like tensors as in the introduction. The decomposition 
$\tn{Ker}\,b = \mathcal W \oplus \Psi(S_0^2(T^*M)) \oplus \Psi(\mbb Rg)$ 
implies that the Weyl tensor vanishes. \end{proof}

\beg{thm} The 4-manifolds $S^2\times \Sigma_g$ with their standard product metric is LCF and hence SD for $g\geq 2$.\end{thm}
\beg{proof} We adapt from \cite{besse}. Suppose we have the constant sectional curvature $+1$ metric on $S^2$ and $-1$ metric on $\Sigma_g$. 
Then we have the following descriptions
$$\tn{Rm}_{S^2}={1\over 2}\,g_{S^2}\varowedge g_{S^2} \hspace{3mm} \tn{and} \hspace{3mm} 
\tn{Rm}_{\Sigma_g}={-1\over 2}\,g_{\Sigma_g}\varowedge g_{\Sigma_g}$$
for the Riemann curvature tensors.
$$\renewcommand{\arraystretch}{2}
\beg{array}{rcl}
\tn{Rm}_{S^2\times \Sigma_g} & = & \tn{Rm}_{S^2} + \tn{Rm}_{\Sigma_g} \\ 
  & = & {1\over 2}\, ( g_{S^2}\varowedge g_{S^2} - g_{\Sigma_g}\varowedge g_{\Sigma_g})\\
    &=&{1\over 2}\, (g_{S^2}+g_{\Sigma_g})\varowedge(g_{S^2}-g_{\Sigma_g})\\
  &=&{1\over 2}\,  g_{S^2 \times \Sigma_g}\varowedge(g_{S^2}-g_{\Sigma_g})\\  
  &=& \Psi({1\over 2}\, (g_{S^2}-g_{\Sigma_g})).  
\end{array}$$
Being in the image of $\Psi$, the Weyl tensor vanishes.
Alternatively, starting with the K\"ahler metrics, one obtains a scalar-flat-K\"ahler (SFK) metric on the product. 
Now apply  Theorem~\ref{sfk=>asd} with both orientations.
\end{proof}

%\newpage

\appendix \section{Appendix}

%Using similar techiques one can also prove the following.
Using Gauss-Bonnet, and signature formula techniques we can also prove the following.

\beg{thm}\label{lemma2KEandLCFchi} If a $4$-manifold admits a K\"ahler-Einstein(KE) metric which is also locally conformally flat(LCF) then its Euler characteristic $\chi=0$.
\end{thm}
\beg{proof} Recall the pointwise identity of Corollary \ref{kahlerWplusandscalarnorms} for K\"ahler metrics    
\beg{equation}  |W_+|^2={s^2 \over 24}. \nonumber \end{equation}
For Einstein metrics, by definition the trace-free Ricci tensor $\stackrel{\circ}{\tn{Ric}}\, =0$.
Plugging these two identities into the Gauss-Bonnet formula (\ref{euler}) we get
\beg{equation}8\pi^2\chi(M)=2\,\|W_+\|^2+\|W_-\|^2. \nonumber \end{equation}
Locally conformally flatness implies $W_\pm=0$, hence $\chi=0$ by above. \end{proof}

As an application we can prove the following non-existence result.

\beg{thm} The product metric on the $4$-manifolds $\Sigma_g\times\Sigma_h$, product of surfaces of genus $g,h\geq 2$ is not a LCF nor SD metric.   
\end{thm}
\beg{proof} First of all, $\Sigma_g\times\Sigma_h$ admits a K\"ahler-Einstein metric. One can see this through different ways. One is to use Aubin/Yau theorem, since $c_1<0$, a surface of general type, there exists a unique KE metric on this complex surface. Another way is to think in terms of product metrics. If you have the hyperbolic $-1$ curvature K\"ahler metrics on both components, the product metric is K\"ahlerian. Besides that, the product of Einstein metrics with common cosmological constant is again Einstein. Combining the two, we obtain a K\"ahler-Einstein metric on the manifold. If the manifold admits a LCF metric as well, then Theorem \ref{lemma2KEandLCFchi} implies that $\chi=0$, however $\chi=(2-2g)(2-2h)$ which is a contradiction. So the product of two hyperbolic metrics on $\Sigma_g\times\Sigma_h$ is not LCF. 

To analyze SD structure we need to find the signature of the manifold. A two dimensional oriented surfaces of genus $g$ always has an orientation reversing diffeomorphism (involution) $R_g$. One can construct this by using a mirror reflection or reflection through a point after arranging the holes symmetrically. Then $(R_g,I)$ is going to be an orientation reversing diffeomorphism of the 4-manifold. Since the Hirzebruch signature is diffeomorphism invariant and changing the orientation changes its sign, we have $\tau=-\tau$ hence $\tau(\Sigma_g\times\Sigma_h)=0$. Alternatively one can compute the intersection matrix as $gh(-H_4)\oplus H_2$ and hence the characteristic polynomial $(\lambda^2-1)^{2gh}(\lambda^2-1)$. (Another approach might be exploiting only the parity of the intersection form and use Rokhlin's theorem of divisibility of the signature by $8$ in the case of even intersection forms to obtain at least some of the cases.) If the product metric is a SD metric on the manifold then it is LCF by Theorem \ref{signature0}, which is already violated. \end{proof}
Existence of LCF metrics on the product 
$\Sigma_g\times\Sigma_h$ of surfaces of genus  $g\geq 2$ and  $h\geq 1$ still remains 
as an open problem.

%\newpage
\vspace{.5cm}

%{\small
%\beg{flushleft}
%\textsc{Tuncel\'i \" Un\' ivers\'ites\'i, Turkia.}\\
%\textit{E-mail address :} \texttt{\textbf{kalafg@gmail.com}}
%\end{flushleft}
%}

{\small
\beg{flushleft}
\textsc{%Mathematics Department, 
Michigan State University, 
East Lansing, MI 48824, USA}\\
\textit{E-mail address :} \texttt{\textbf{kalafat@math.msu.edu}}
\end{flushleft}
}

%\vspace{1cm}

\end{document}